\documentclass[preprint,3p,12pt,pdf]{elsarticle}

\usepackage{lineno,hyperref}
\modulolinenumbers[5]

\usepackage{hyperref}

\usepackage{epstopdf}

\usepackage{tikz}
\usetikzlibrary{arrows}

\usepackage{pst-node}
\usepackage{tikz-cd}

\usepackage[shellescape]{gmp}

\usepackage{mathrsfs}
\usepackage{amssymb}
\usepackage{amsfonts}
\usepackage{latexsym}
\usepackage{mathtools}
\usepackage{xcolor}
\usepackage{wrapfig}
\usepackage{floatflt}

\usepackage{mathtools}
\usepackage{extarrows}

\usepackage{graphicx}

\usepackage{subcaption}

\def\lb{\label}

\newcommand{\er}[1]{\textrm{(\ref{#1})}}

\def\a{\alpha}

   \def\vt{\vartheta}    \def\vp{\varphi}    

\def\Z{{\mathbb Z}}    \def\R{{\mathbb R}}   \def\C{{\mathbb C}}    
    \def\N{{\mathbb N}}   

\def\lt{\biggl}                  \def\rt{\biggr}
\def\ol{\overline}               \def\wt{\widetilde}

\let\ge\geqslant                 \let\le\leqslant

\def\iy{\infty}

\def\pa{\partial}                
                 \def\ev{\equiv}
        
\def\el2{\ell^{\,2}}             \def\1{1\!\!1}

\def\Im{\mathop{\mathrm{Im}}\nolimits}

\def\Re{\mathop{\mathrm{Re}}\nolimits}

\let\ge\geqslant
\let\le\leqslant

\newcommand{\ca}{\begin{cases}}
\newcommand{\ac}{\end{cases}}
\newcommand{\ma}{\begin{pmatrix}}
\newcommand{\am}{\end{pmatrix}}
\def\eq{\begin{equation}}
\def\qe{\end{equation}}
\def\[{\begin{equation}}
\def\]{\end{equation}}

\begin{document}

\begin{frontmatter}

\title{Approximation of the number of descendants in branching processes}

\date{\today}

\author
{Anton A. Kutsenko}

\address{Jacobs University, 28759 Bremen, Germany; email: akucenko@gmail.com}

\begin{abstract}
We discuss approximations of the relative limit densities of descendants in Galton--Watson processes that follow from the Karlin--McGregor near-constancy phenomena. These approximations are based on the fast exponentially decaying Fourier coefficients of Karlin--McGregor functions and the binomial coefficients. The approximations are sufficiently simple and show good agreement between approximate and exact values, which is demonstrated by several numerical examples.
\end{abstract}



\end{frontmatter}


{\section{Introduction}\lb{sec1}}

Near-constancy phenomena indicates that the tail of limit distributions of simple branching processes has a power-law asymptotic perturbed by extremely small oscillations. This phenomena has a long history starting from the works \cite{KM1,KM2} and then continues to be actively studied in \cite{D,BB,FW}. It should also be noted interesting and important physical applications, see \cite{DIL,CG}. Most of the research is devoted to the {\it martingale limit}, see, e.g. \cite{S}, for branching processes, which represents a continuous function of a real argument. For this function, the power-law asymptotic multiplied by some multiplicatively periodic factor is proven. In our research, we focus on Taylor coefficients of another limit distribution, which satisfies some Schr\"oder functional equation. The analysis of the asymptotic of the power-series coefficients of solutions of functional equations is difficult and may strongly depends on the functional equation itself, see \cite{O,B}. For the sufficiently general type of Schr\"oder functional equations, we provide an approximation formula for the Taylor coefficients of the solutions of these equations. The approximation formula is based on the special binomial coefficients, from asymptotic of which one may extract the leading power-law term and periodic factors with small amplitudes. 

The Galton--Watson process is defined by
\[\lb{001}
 X_{t+1}=\sum_{j=1}^{X_t}\xi_{j,t},\ \ \ X_0=1,\ \ \ t\in\N\cup\{0\},
\]
where all $\xi_{j,t}$ are independent and identically-distributed natural number-valued random variables with the probability generating function
\[\lb{002}
 P(z):=\mathbb{E}z^{\xi}=p_0+p_1z+p_2z^2+p_3z^3+....
\]
It is well known, see \cite{H}, that 
\[\lb{003}
 \mathbb{E}z^{X_t}=\underbrace{P\circ...\circ P}_{t}(z),\ \ \ t\in\N.
\]
Let us assume $p_0=0$ and $0<p_1<1$. In this case we can define
\[\lb{004}
 \Phi(z)=\lim_{t\to\iy}p_1^{-t}\underbrace{P\circ...\circ P}_{t}(z),
\]
which is analytic at least for $|z|<1$. In fact, $\Phi(z)$ is analytic inside a component of the filled Julia set for $P(z)$ which contains the ball $B_1:=\{z\in\C:\ |z|<1\}$, since $|P(z)|\le|z|$ for $z\in B_1$. The function $\Phi$ satisfies the Scr\"oder-type functional equation
\[\lb{005}
 \Phi(P(z))=p_1\Phi(z),\ \ \ \Phi(0)=0,\ \ \ \Phi'(0)=1.
\] 
Due to \er{002}-\er{004}, the Taylor coefficients $\varphi_n$ of
\[\lb{006}
 \Phi(z)=z+\varphi_2z^2+\varphi_3z^3+...
\]
describe the relative limit densities of the number of descendants 
\[\lb{007}
 \varphi_n=\lim_{t\to+\iy}\frac{\mathbb{P}\{Z_t=n\}}{\mathbb{P}\{Z_t=1\}}.
\]
Our goal is to find simple approximations of $\varphi_n$. Let us assume that the first moment
\[\lb{008}
 E:=P'(1)=\sum_{j=1}^{+\iy}jp_j<+\iy
\]
is finite. Note that $E>1$, since $\sum_{j=1}^{+\iy}p_j=1$ and $0<p_1<1$. Thus, we can define
\[\lb{009}
 \Psi(z)=\lim_{t\to\iy}E^t(1-\underbrace{P^{-1}\circ...\circ P^{-1}}_{t}(z)),
\]
which is analytic in some neighborhood of $z=1$. In fact, it is not difficult to show that $\Psi$ is analytic on $(0,1)$, including some its neighborhood, since $P$ is strictly monotonic on this interval and $P(0)=0$, $P(1)=1$. The function $\Psi$ also satisfies the Schr\"oder-type functional equation
\[\lb{010}
 \Psi(P(z))=E\Psi(z),\ \ \ \Psi(1)=0,\ \ \ \Psi'(1)=-1.
\]
Everything is ready to define 
\[\lb{011}
 \Theta(z)=\Phi(z)\Psi(z)^{\frac{-\ln p_1}{\ln E}},
\]
which is analytic on $(0,1)$. Moreover, it satisfies the functional equation
\[\lb{012}
 \Theta(P(z))=\Theta(z),
\]
see \er{005}, \er{010}, and \er{011}. The ``polynomially" periodic function $\Theta$ can be considered as some precursor of the multiplicatively periodic Karlin--McGregor function $K(z)=\Theta(\Psi^{-1}(z))$, 
see \er{010} and \er{012}. The function $\Psi^{-1}$ can be defined as
\[\lb{013}
 \Pi(z):=\Psi^{-1}(z)=\lim_{t\to+\iy}\underbrace{P\circ...\circ P}_{t}(1-E^{-t}z),
\]
see \er{009}. The function $\Pi(z)$ is an entire function which satisfies the Poincar\'e-type functional equation
\[\lb{014}
 P(\Pi(z))=\Pi(Ez),\ \ \ \Pi(0)=1,\ \ \ \Pi'(0)=-1.
\]
The next step is to introduce
\[\lb{015}
 K^*(z)=K(E^z)=\Theta(\Pi(E^z))=\Phi(\Pi(E^z))p_1^{-z},
\] 
which is an $1$-periodic function
defined on some strip $\{z:\ |\Im z|<\vt^*\}$ with $\vt^*>0$. Thus, we have the Fourier series
\[\lb{016}
 K^*(z)=\sum_{m=-\iy}^{+\iy}\theta_me^{2\pi \mathbf{i} mz},\ \ \ |\Im z|<\vt^*
\] 
with exponentially rapidly decreasing coefficients $\theta_m=\ol{\theta_{-m}}$. We denote $\mathbf{i}=\sqrt{-1}$. Using the definitions $\Psi^{-1}=\Pi$ and \er{015}, along with \er{016} we obtain
\[\lb{017}
 \Theta(z)=\sum_{m=-\iy}^{+\iy}\theta_m\Psi(z)^{\frac{2\pi \mathbf{i} m}{\ln E}},
\] 
which with \er{011} leads to
\[\lb{018}
 \Phi(z)=\sum_{m=-\iy}^{+\iy}\theta_m\Psi(z)^{\frac{2\pi \mathbf{i} m+\ln p_1}{\ln E}}.
\]
The main idea is that the amplitudes of oscillations $\theta_m$ are usually extremely small. Thus, we can take few terms to obtain a good approximation
\[\lb{019}
 \Phi(z)\approx\sum_{m=-M}^{M}\theta_m\Psi(z)^{\frac{2\pi \mathbf{i} m+\ln p_1}{\ln E}}.
\]
After that, we substitute the expansion of $\Psi(z)$ near $z=1$, see \er{009} and \er{010},
\[\lb{020}
 \Psi(z)=1-z+O((1-z)^2)\ \ {\rm or}\ \ \Psi(z)=1-z+\frac{P''(1)}{2(E^2-E)}(1-z)^2+...+O((1-z)^R)
\]
into \er{019} to obtain
\[\lb{021}
 \Phi(z)\approx\sum_{m=-M}^{M}\theta_m(1-z+...)^{\frac{2\pi \mathbf{i} m+\ln p_1}{\ln E}}.
\]
In comparison to the computations of $\varphi_n$, the Taylor coefficients of RHS in \er{021} can be computed easily. For example, if we take the first approximation in \er{020} only then \er{021} becomes
\[\lb{022}
 \Phi(z)\approx\sum_{m=-M}^{M}\theta_m(1-z)^{\frac{2\pi \mathbf{i} m+\ln p_1}{\ln E}}=\sum_{n=0}^{+\iy}z^n\sum_{m=-M}^M(-1)^n\binom{\frac{2\pi \mathbf{i} m+\ln p_1}{\ln E}}{n}\theta_m,
\]
which gives an approximation
\[\lb{023}
 \varphi_n\approx\tilde\vp_n^M:=(-1)^n\sum_{m=-M}^M\binom{\frac{2\pi \mathbf{i} m+\ln p_1}{\ln E}}{n}\theta_m,
\]
where the generalized binomial coefficients 
\[\lb{024}
 \binom{k}{n}=\frac{k(k-1)(k-2)...(k-n+1)}{n!},\ \ \ n\ge1 
\]
are Taylor coefficients of $(1+z)^k$ for arbitrary $k$. Using the asymptotics
\[\lb{025}
 \binom{k}{n}\approx\frac{(-1)^n}{\Gamma(-k)n^{k+1}},\ \ \ n\to+\iy,
\]
and the Stirling approximation 
\begin{multline}\lb{026}
 \Gamma(-k)\approx\sqrt{2\pi}\exp\lt(\lt(-k-\frac12\rt)\ln(-k) +k\rt)=\sqrt{2\pi}\exp\lt(\lt(\frac{-2\pi \mathbf{i} m-\ln p_1}{\ln E}-\frac12\rt)\cdot\\
 \cdot\ln\frac{-2\pi \mathbf{i} m-\ln p_1}{\ln E}+\frac{2\pi \mathbf{i} m+\ln p_1}{\ln E}\rt)\approx C_{m}e^{-\frac{\pi^2|m|}{\ln E}}|m|^{-\frac{\ln p_1}{\ln E}-\frac12}
\end{multline}
for some bounded constants $C_{m}$, and $k=\frac{2\pi \mathbf{i} m+\ln p_1}{\ln E}$, see \er{023}, with $m\to\iy$, we obtain the approximation
\[\lb{027}
 \binom{\frac{2\pi \mathbf{i} m+\ln p_1}{\ln E}}{n}\theta_m\approx  C_m\theta_me^{\frac{\pi^2|m|}{\ln E}}|m|^{\frac{\ln p_1}{\ln E}+\frac12},
\] 
for
large $n$ and $m$. While this approximation is very rough and we do not use it directly, it allows us to estimate the possible order of growths of terms in \er{023}. Below, we discus some other details related to \er{027}.

{\bf Remark on exponential decaying $\theta_m$.} It is useful to note the following remark. If $z=x+\mathbf{i}y$ for $x\to-\iy$ and $y\in\R$ then $\Pi(E^z)=1-E^xe^{\mathbf{i}y\ln E}+o(E^x)$, since $\Pi(0)=1$ and $\Pi'(0)=-1$. Thus, 
\[\lb{029}
 K^*(z)=\Phi(1-E^xe^{\mathbf{i}y\ln E}+o(E^x))p_1^{-x-\mathbf{i}y}. 
\] 
As mentioned above, the function $\Phi$ is defined on the filled Julia set for the polynomial $P$ which contains the unit ball. Hence, the maximal $y$ for which $K^*(z)$ is defined is at least $\frac{\pi}{2\ln E}$ for $x\to-\iy$, see \er{029}, since the unit circle intersects real axis at the angle $\frac{\pi}2$. Thus, the maximal width of the strip where $1$-periodic function $K^*$ is defined satisfies $\vt^*\ge\frac{\pi}{2\ln E}$. Moreover, at the edge of the smaller strip the function $K^*$ is bounded by some value $C_P$, since the fractal open filled Julia set for the polynomial $P$ contains a segment of the line $\{1+\mathbf{i}y,\ y\in\R\}$ located near $z=1$, except this unique point $z=1$. For instance, this can be verified by observing a polynomial $P$ iteration of the points of this segment. For small $y$, each next iteration decreases the real part of the points 
\[\lb{030}
 P(1+\mathbf{i}y)=1+E\mathbf{i}y-\frac{P''(1)}{2}y^2+O(y^3),
\]
and at some step the points are inside the unit ball. The modified Karlin-McGregor function
\[\lb{031}
 K^*\lt(\frac{\ln z}{2\pi \mathbf{i}}\rt)=\sum_{m=-\iy}^{+\iy}\theta_mz^m
\]
is analytic on the ring $\{z:\ e^{-2\pi\vt^*}<|z|<e^{2\pi\vt^*}\}$. Using $\vt^*\ge\frac{\pi}{2\ln E}$ and applying Cauchy estimates for this smaller ring we obtain
\[\lb{032}
 |\theta_{m}|\le C_Pe^{-\frac{|m|\pi^2}{\ln E}},\ \ \ m\in\Z,
\]
which shows rapid exponential decaying of $\theta_m$. Unfortunately, the binomial factor in \er{023} may compensate this exponential decaying, see \er{027}. However, the factor $|m|^{\frac{\ln p_1}{\ln E}+\frac12}$ in \er{027} and, especially, the exact value of $\vt^*$, which often greater than $\frac{\pi}{2\ln E}$, allow us to hope for a very rapid decrease of terms in \er{023}. (Note that $\er{023}$ converges anyway, see \er{106} below.) As an example, we draw the Julia set for the polynomial $P(z)=0.1z+0.5z^2+0.4z^3$, see Fig. \ref{fig1}, where the so-called ``critical angle" is strictly less than $\pi$, which leads to $\vt^*>\frac{\pi}{2\ln E}$, see explanations above \er{030}. The ``critical angle" is the minimal angle in the neighborhood of $1$, exterior of which asymptotically belongs to the filled Julia set when the radius tends to $0$. 

For practical implementations, it is useful to compute Fourier coefficients of $K^*(x-\frac{\mathbf{i}\pi}{2\ln E})$ instead of $K^*(x)$, because they contain already the exponential factor $e^{\frac{m\pi^2}{\ln E}}$. After that, $\frac{m\pi^2}{\ln E}$ should be subtracted from $\ln\binom{\frac{2\pi \mathbf{i} m+\ln p_1}{\ln E}}{n}$ before taking exponent. This procedure leads to stable and accurate computations of the terms in \er{023} for positive $m$. It is not necessary to compute the terms with negative $m$, since they are conjugations of ``positive" ones. We actively use this trick for the computations in the examples Section.    

\begin{figure}[h]
	\center{\includegraphics[width=0.9\linewidth]{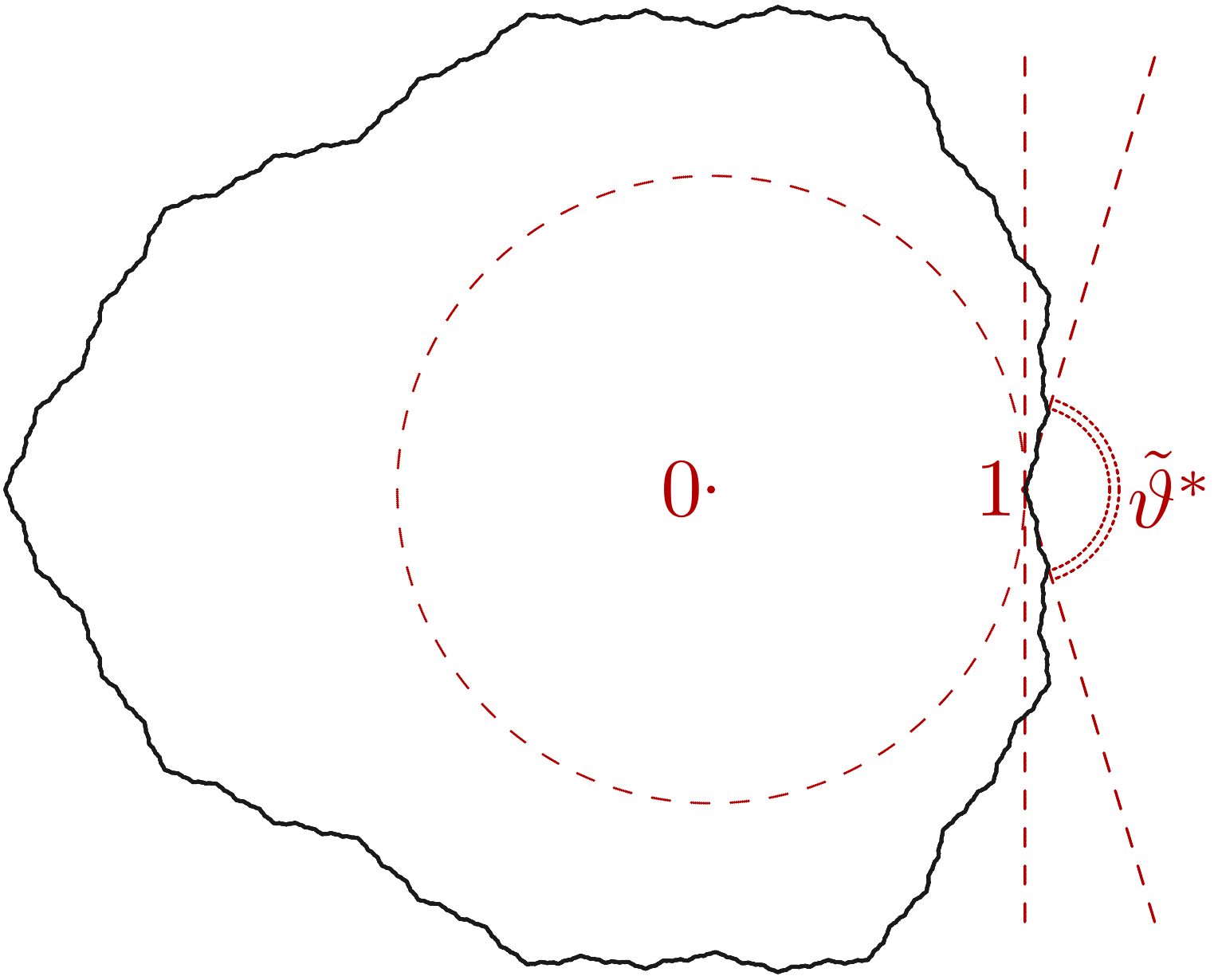}}
	\caption{Julia set for the polynomial $P(z)=0.1z+0.5z^2+0.4z^3$. It is the boundary of the open filled Julia set, which is the domain of definition for the analytic function $\Phi$, see \er{004} and \er{005}. The unit ball belongs to the filled Julia set. The critical angle $\tilde\vt^*=2(\pi-\vt^*\ln E)$ is computed approximately.}\lb{fig1}
\end{figure}

{\bf Remark on computation of $\theta_m$} Approximations \er{019}-\er{023} are based on the computation of $\theta_m$, $m\in\Z$, or also $e^{\frac{\pi^2m}{\ln E}}\theta_m$, $m\in\N$. Thus we should compute $K^*(z)$, since it is periodic and
\[\lb{100}
 \theta_m=\int_0^1K^*(x)e^{-2\pi i nx}dx\ \ {\rm or}\ \ e^{\frac{\pi^2m}{\ln E}}\theta_m=\int_0^1K^*(x-\frac{\mathbf{i}\pi}{2\ln E})e^{-2\pi i mx}dx.
\]
In the end, the calculation of 
\[\lb{101}
 K^*(z)=\Phi(\Pi(E^z)){p_1}^{-z}
\]
is reduced to the calculation of $\Phi$ and $\Pi$, see \er{015}. Let us assume that $P$, see \er{002}, is a polynomial
\[\lb{102}
 P(z)=p_1z+p_2z^2+...+p_Nz^N,\ \ \ 0<p_1<1,
\]
which is a most practical case. By induction, let us define
\begin{multline}\lb{103}
 \Phi_0(z)=z,\ \ \ \Phi_{t+1}(z)=p_1^{-t-1}\underbrace{P\circ...\circ P}_{t+1}(z)=
 p_1^{-t}p_1^{-1}P(\underbrace{P\circ...\circ P}_{t}(z))=\\
 \Phi_t(z)+\sum_{j=2}^Np_jp_1^{(j-1)t-1}\Phi_t(z)^j,
\end{multline}
see \er{102}. The recurrent formula \er{103} provides an exponentially fast convergence $\Phi_t(z)\to\Phi(z)$, since $p_1^{(j-1)t-1}\to0$ for $t\to+\iy$ exponentially fast when $j\ge2$.

To compute $\Pi(z)$, we rewrite \er{102} as
\[\lb{104}
 P(z)=1-E(1-z)+q_2(1-z)^2+...+q_N(1-z)^N,
\]
see \er{008}. By induction, let us define
\begin{multline}\lb{105}
 \Pi_0(z)=z,\ \ \ \Pi_{t+1}(z)=\underbrace{P\circ...\circ P}_{t+1}(1-E^{-t-1}z)=\underbrace{P\circ...\circ P}_{t}(P(1-E^{-t-1}z))=\\
 \underbrace{P\circ...\circ P}_{t}(1-E^{-t}(z-\sum_{j=2}^Nq_jE^{-(j-1)t-j}z^j))=\Pi_t(z-\sum_{j=2}^Nq_jE^{-(j-1)t-j}z^j)).
\end{multline}
The type of the recurrent formula \er{105} differs from that of \er{103}, but \er{105} also provides the exponentially fast convergence $\Pi_t(z)\to\Pi(z)$ for $t\to+\iy$, since $E^{-(j-1)t-j}\to0$ for $t\to+\iy$ exponentially fast when $j\ge2$. In the examples below we use numerical procedures based on \er{103} and \er{105}. 

{\bf Remark on the approximations for $M\to+\iy$.} For fixed $n$, the series \er{023} and \er{203} converges when $M\to+\iy$, since $\theta_m$ decay exponentially fast, while the binomial coefficients are polynomials of $m$, see \er{024}. Using \er{016} and \er{023}, one can express the best in the class approximation $\wt\vp_n^{\iy}$ through the derivatives of $K^*$, namely
\[\lb{106}
 \tilde\vp_n^{\iy}=(-1)^n\binom{\frac{\frac{\pa}{\pa z}+\ln p_1}{\ln E}}{n}K^*(z)|_{z=0},
\]
since the derivative of $K^*$ leads to the multiplication of each coefficient $\theta_m$ by $2\pi\mathbf{i}m$. 

{\bf Remark on $\vp_n$.} It is necessary to take into account the well known fact that, in some cases, we should be careful using the approximations. For example, if all $p_{2n}=0$ then all $\vp_{2n}=0$ as well. In such cases, the approximation procedure \er{019}-\er{023} should be modified. The modification is not very difficult, but, for simplicity, let us assume that $P(z)$ is a polynomial with all non-zero coefficients except $p_0=0$.

We now move on to examples in Section \ref{sec2}, the general expansion will be considered in Section \ref{sec3}, the conclusion will be given in Section \ref{sec4}. 

{\section{Examples}\lb{sec2}}

If $\frac{\ln p_1}{\ln E}>-2$ then $M=0$ or $M=1$ in \er{023} already give good approximation of $\vp_n$. The more interesting case is when $p_1$ is small. Let us assume that
\[\lb{200}
 -3<\frac{\ln p_1}{\ln E}\le-2.
\]
In this case, we take two terms in the Taylor expansion for $\Psi(z)$, since $(1-z)^{\frac{2\pi\mathbf{i}m+\ln p_1}{\ln E}}$ and $(1-z)^{1+\frac{2\pi\mathbf{i}m+\ln p_1}{\ln E}}$ give power-series with non-attenuating coefficients, while the coefficients of $(1-z)^{2+\frac{2\pi\mathbf{i}m+\ln p_1}{\ln E}}$ already tend to $0$. Thus, we have
\begin{multline}\lb{201}
 \Psi(z)^{\frac{2\pi\mathbf{i}m+\ln p_1}{\ln E}}=\lt(1-z+\frac{P''(1)}{2(E^2-E)}(1-z)^2+O((1-z)^3)\rt)^{\frac{2\pi\mathbf{i}m+\ln p_1}{\ln E}}=(1-z)^{\frac{2\pi\mathbf{i}m+\ln p_1}{\ln E}}\cdot\\
 \cdot\lt(1+\frac{P''(1)}{2(E^2-E)}(1-z)+O((1-z)^2)\rt)^{\frac{2\pi\mathbf{i}m+\ln p_1}{\ln E}}=(1-z)^{\frac{2\pi\mathbf{i}m+\ln p_1}{\ln E}}+\\
 +\frac{P''(1)(2\pi\mathbf{i}m+\ln p_1)}{2(E^2-E)\ln E}(1-z)^{1+\frac{2\pi\mathbf{i}m+\ln p_1}{\ln E}}+O((1-z)^{2+\frac{2\pi\mathbf{i}m+\ln p_1}{\ln E}}),
\end{multline}
which leads to the approximation
\[\lb{202}
 \Phi(z)\approx\sum_{m=-M}^M\lt((1-z)^{\frac{2\pi\mathbf{i}m+\ln p_1}{\ln E}}
 +\frac{P''(1)(2\pi\mathbf{i}m+\ln p_1)}{2(E^2-E)\ln E}(1-z)^{1+\frac{2\pi\mathbf{i}m+\ln p_1}{\ln E}}\rt)\theta_m,
\]
see \er{019}, which, in turn, leads to
\[\lb{203}
 \vp_n\approx\vp_n^M:=(-1)^n\sum_{m=-M}^M\lt(\binom{\frac{2\pi\mathbf{i}m+\ln p_1}{\ln E}}{n}
 +\frac{P''(1)(2\pi\mathbf{i}m+\ln p_1)}{2(E^2-E)\ln E}\binom{1+\frac{2\pi\mathbf{i}m+\ln p_1}{\ln E}}{n}\rt)\theta_m.
\]
It is useful to compare \er{203} with \er{023}, where the new correction terms appear. These new terms are not leading, but they improve the approximation. 

It is possible to extract the conceptual behavior - the power factor and the oscillations - of $\vp_n^M$ from \er{203}. It is enough to take two terms in the asymptotic of binomial coefficients
\[\lb{203n0}
 \binom{k}{n}=\frac{(-1)^n}{\Gamma(-k)n^{k+1}}+\frac{(-1)^nk(k+1)}{2\Gamma(-k)n^{k+2}}+...,
\]
since the next terms already tend to zero for $n\to\iy$ and $\Re k>-3$. Using \er{203n0} in \er{203} along with $\Gamma(k+1)=k\Gamma(k)$, we obtain
\begin{multline}\lb{203n1}
 \vp_n\approx\vp_n^{\Gamma,M}:=n^{\frac{-\ln p_1}{\ln E}-1}\sum_{m=-M}^M\frac{\theta_me^{2\pi\mathbf{i}m\frac{-\ln n}{\ln E}}}{\Gamma(-\frac{2\pi\mathbf{i}m+\ln p_1}{\ln E})}
 +\\ \frac{n^{\frac{-\ln p_1}{\ln E}-2}(P''(1)+E-E^2)}{2(E^2-E)\ln E}\sum_{m=-M}^M\frac{(2\pi\mathbf{i}m+\ln p_1)\theta_me^{2\pi\mathbf{i}m\frac{-\ln n}{\ln E}}}{\Gamma(-\frac{2\pi\mathbf{i}m+\ln p_1}{\ln E}-1)}.
\end{multline}
In fact, \er{203n1} consists of the terms $n^{\a}{\rm Periodic}(\log_E n)$, which is very useful for understanding the conceptual behavior of $\vp_n$.

Let us discuss how to compute exact values $\vp_n$ for the cubic polynomial $P(z)=p_1z+p_2z^2+p_3z^3$. Substituting \er{006} into \er{005} we obtain
\[\lb{204}
 \sum_{n=1}^{+\iy}\vp_n(p_1z+p_2z^2+p_3z^3)^n=p_1\sum_{n=1}^{+\iy}\vp_nz^n,\ \ \ \vp_1=1,
\]
which can be rewritten in the form
\[\lb{205}
  \sum_{n=1}^{+\iy}\vp_n\sum_{i+j+k=n}\frac{n!p_1^ip_2^jp_3^kz^{i+2j+3k}}{i!j!k!}=p_1\sum_{n=1}^{+\iy}\vp_nz^n,\ \ \ \vp_1=1,
\]
which, in turn, leads to
\[\lb{206}
 \sum_{n=1}^{+\iy}\vp_n\sum_{r=n}^{3n}c_{n,r}z^r=p_1\sum_{n=1}^{+\iy}\vp_nz^n,\ \ \ \vp_1=1,\ \ \ c_{n,r}=\sum_{k=0}^{\frac{r-n}2}\frac{n!p_1^{2n+k-r}p_2^{r-n-2k}p_3^k}{(2n+k-r)!(r-n-2k)!k!},
\]
which, finally, gives the recurrent formula
\[\lb{207}
 \vp_n=\frac1{p_1-p_1^n}\sum_{j=1}^{\frac{2n}3}c_{n-j,n}\vp_{n-j},\ \ \ \vp_1=1,
\]
that allows us to compute all $\vp_n$ recurrently. Yes, even for cubic polynomials, the computation of exact values $\vp_n$ is much harder than the computation of their approximations, as it is seen from the comparison of \er{206} and \er{207} with \er{203} and \er{203n1}.

\subsection{Example $p_1=0.2$, $p_2=0.6$, $p_3=0.2$.} We plot exact values \er{207} and their approximations \er{203} and \er{023} for the cubic polynomial $P(z)=0.2z+0.6z^2+0.2z^3$, see Figs. \ref{fig2} and \ref{fig3}. The simple approximations $\vp_n^0$ and $\tilde\vp_n^0$ give already good results, see Fig. \ref{fig2}.(a). They show the dominant power-law growth $\frac{\theta_0n^{-1-\ln p_1/\ln E}}{\Gamma(-\ln p_1/\ln E)}$, see \er{023} and \er{025}. However, $\vp_n^1$ and, especially, $\vp_n^2$ are much better than $\vp_n^0$ and, of course, than $\tilde\vp_n^0$. Namely, comparing Fig. \ref{fig3}.(a) and (b), it is seen that the additional terms in \er{023} play a very important role. The accurate $\Gamma$-approximations are also very good, see Fig. \ref{fig3gamma}. Let us provide some important characteristics for this example
$$
 \frac{-\ln p_1}{\ln E}\approx2.32,\ \ \ \theta_0\approx1.94,\ \ e^{\frac{\pi^2}{\ln E}}\theta_1\approx0.16-1.3\mathbf{i},\ \ e^{\frac{2\pi^2}{\ln E}}\theta_2\approx-0.0036+0.04\mathbf{i},\ \ e^{\frac{\pi^2}{\ln E}}\approx1.5\cdot10^6.
$$

\begin{figure}[h]
    \centering
    \begin{subfigure}[b]{0.45\textwidth}
        \includegraphics[width=\textwidth]{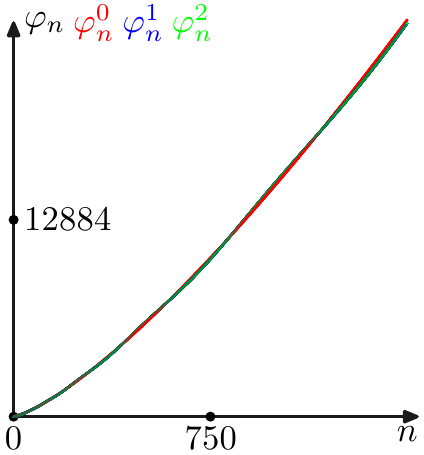}
        \caption{}
    \end{subfigure}
    \begin{subfigure}[b]{0.45\textwidth}
        \includegraphics[width=\textwidth]{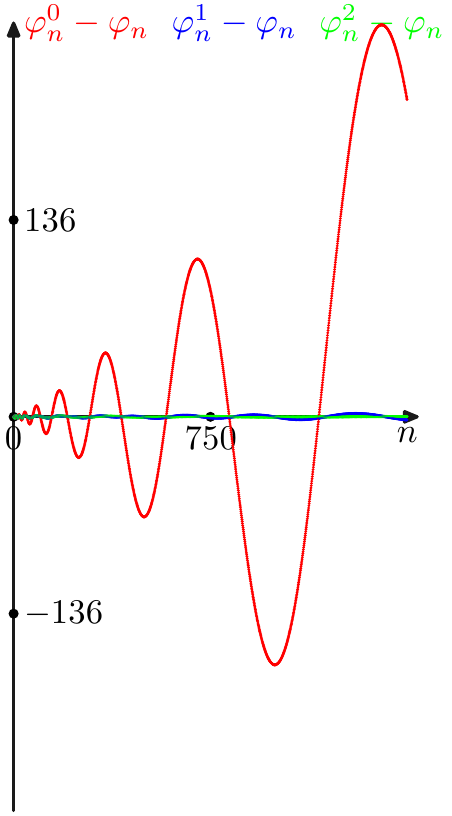}
        \caption{}
    \end{subfigure}
    \caption{For the polynomial $P(z)=0.2z+0.6z^2+0.2z^3$, Tailor coefficients $\vp_n$ (black curve) and their approximations $\vp_n^M$ (red, blue, and green curves) for $M=0,1,2$ are plotted in (a), the differences $\vp_n^M-\vp_n$ (red, blue, and green curves) are plotted in (b).}\lb{fig2}
\end{figure}

\begin{figure}[h]
    \centering
    \begin{subfigure}[b]{0.45\textwidth}
        \includegraphics[width=\textwidth]{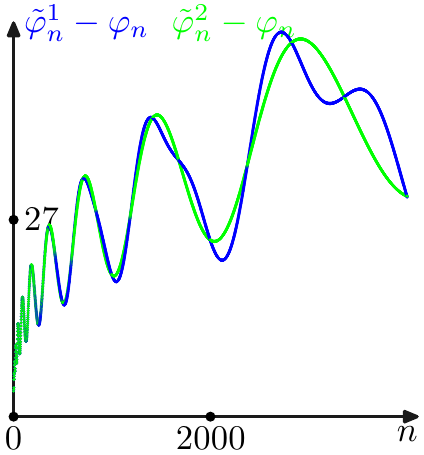}
        \caption{}
    \end{subfigure}
    \begin{subfigure}[b]{0.45\textwidth}
        \includegraphics[width=\textwidth]{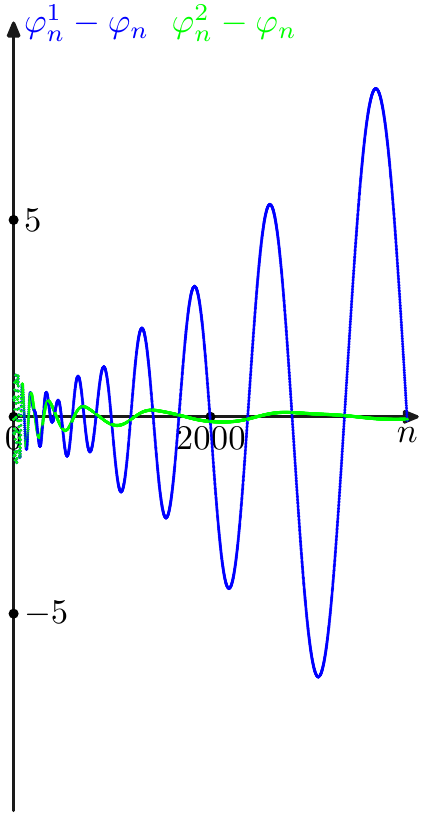}
        \caption{}
    \end{subfigure}
    \caption{For the polynomial $P(z)=0.2z+0.6z^2+0.2z^3$, the differences $\tilde\vp_n^M-\vp_n$ and $\vp_n^M-\vp_n$ (blue and green curves) are plotted in (a) and (b) respectively.}\lb{fig3}
\end{figure}

\begin{figure}[h]
    \centering
    \begin{subfigure}[b]{0.45\textwidth}
        \includegraphics[width=\textwidth]{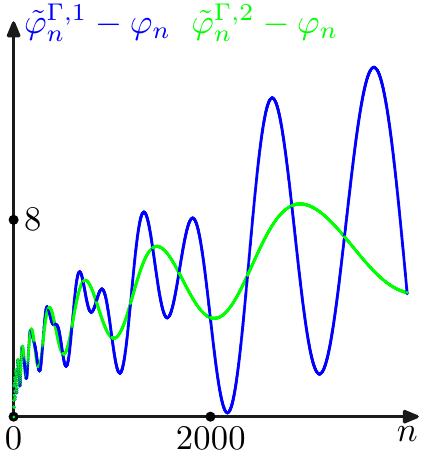}
        \caption{}
    \end{subfigure}
    \begin{subfigure}[b]{0.45\textwidth}
        \includegraphics[width=\textwidth]{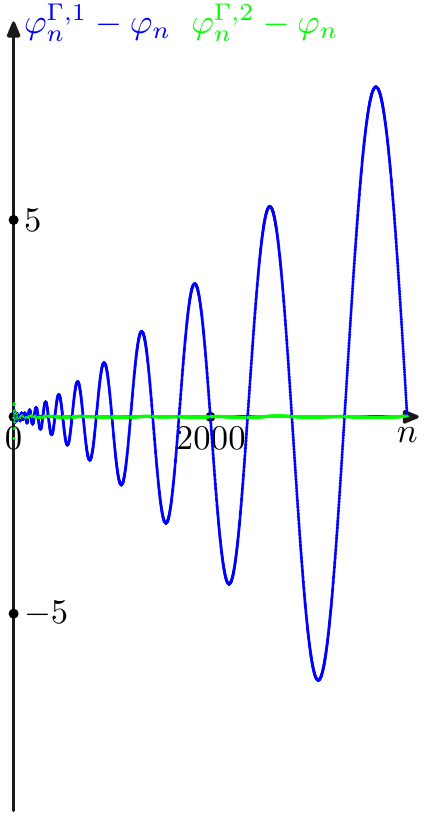}
        \caption{}
    \end{subfigure}
    \caption{For the polynomial $P(z)=0.2z+0.6z^2+0.2z^3$, the differences $\tilde\vp_n^{\Gamma,M}-\vp_n$ and $\vp_n^{\Gamma,M}-\vp_n$ (blue and green curves) are plotted in (a) and (b) respectively. The $\Gamma$-approximations $\tilde\vp_n^{\Gamma,M}$ and $\vp_n^{\Gamma,M}$ are given by \er{023} and \er{203n0} respectively, where, in the first case, the binomial coefficients are replaced by their approximations \er{024}.}\lb{fig3gamma}
\end{figure}

\subsection{Example $p_1=0.1$, $p_2=0.5$, $p_3=0.4$.} As another example, we compute $\vp_n$ and their approximations for the polynomial $P(z)=0.1z+0.5z^2+0.4z^3$, see Figs. \ref{fig4} and \ref{fig5}. Some important characteristics for this example are
$$
 \frac{-\ln p_1}{\ln E}\approx2.76,\ \ \ \theta_0\approx2.887,\ \ e^{\frac{\pi^2}{\ln E}}\theta_1\approx9.94-0.62\mathbf{i},\ \ e^{\frac{2\pi^2}{\ln E}}\theta_2\approx1.21+1.66\mathbf{i},\ \ e^{\frac{\pi^2}{\ln E}}\approx1.4\cdot10^6.
$$

\begin{figure}[h]
    \centering
    \begin{subfigure}[b]{0.45\textwidth}
        \includegraphics[width=\textwidth]{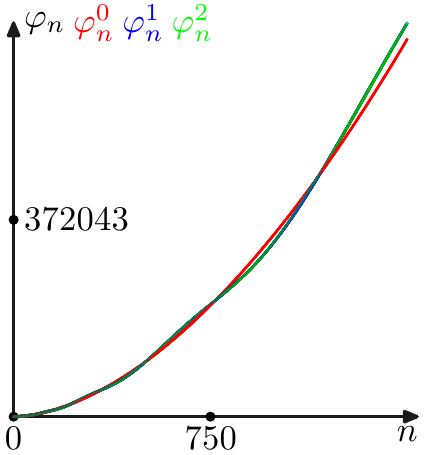}
        \caption{}
    \end{subfigure}
    \begin{subfigure}[b]{0.45\textwidth}
        \includegraphics[width=\textwidth]{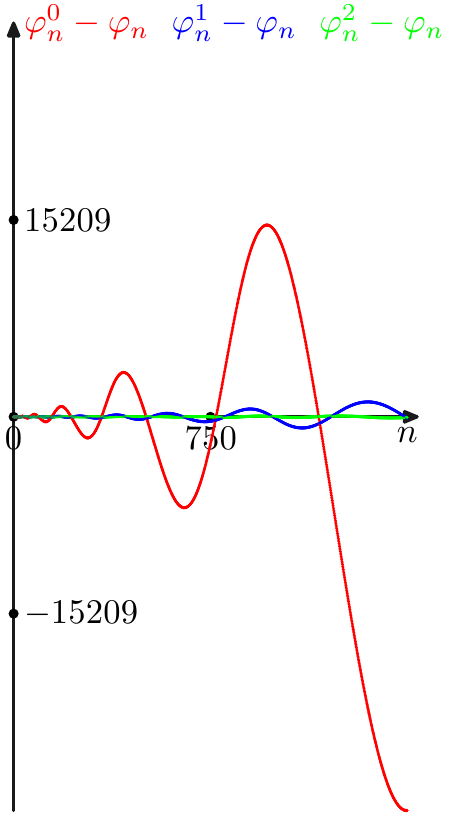}
        \caption{}
    \end{subfigure}
    \caption{For the polynomial $P(z)=0.1z+0.5z^2+0.4z^3$, Tailor coefficients $\vp_n$ (black curve) and their approximations $\vp_n^M$ (red, blue, and green curves) for $M=0,1,2$ are plotted in (a), the differences $\vp_n^M-\vp_n$ (red, blue, and green curves) are plotted in (b).}\lb{fig4}
\end{figure}

\begin{figure}[h]
    \centering
    \begin{subfigure}[b]{0.45\textwidth}
        \includegraphics[width=\textwidth]{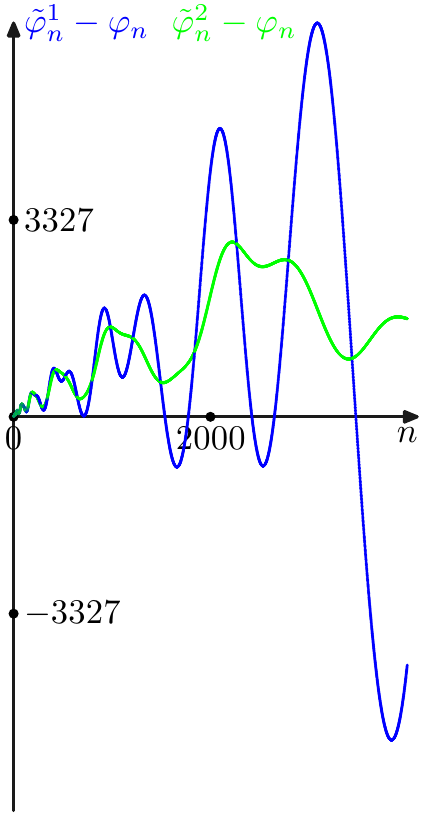}
        \caption{}
    \end{subfigure}
    \begin{subfigure}[b]{0.45\textwidth}
        \includegraphics[width=\textwidth]{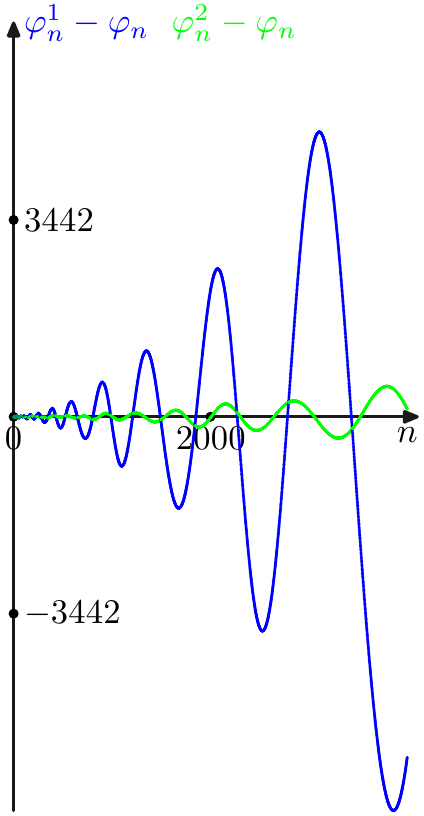}
        \caption{}
    \end{subfigure}
    \caption{For the polynomial $P(z)=0.1z+0.5z^2+0.4z^3$, the differences $\tilde\vp_n^M-\vp_n$ and $\vp_n^M-\vp_n$ (blue and green curves) are plotted in (a) and (b) respectively.}\lb{fig5}
\end{figure}

\subsection{Example $p_1=0.25$, $p_2=0.5$, $p_3=0.25$.} In fact, I made a lot of tests, not only for cubic polynomials. All of them confirm the main observations obtained in the previous examples. To finish this Section, I made one example which complements some particular results from \cite{K1} (see Figure 3 there).  Namely,  the approximations $\vp_n^1$ and $\vp_n^2$ improve significantly $\vp_n^0$, see Fig. \ref{fig6}. This example is my main personal motivation for this work. Some important characteristics for this example are
$$
 \frac{-\ln p_1}{\ln E}=2,\ \ \ \theta_0\approx1.46,\ \ e^{\frac{\pi^2}{\ln E}}\theta_1\approx0.12+0.01\mathbf{i},\ \ e^{\frac{2\pi^2}{\ln E}}\theta_2\approx0.002+0.0001\mathbf{i},\ \ e^{\frac{\pi^2}{\ln E}}\approx1.5\cdot10^6.
$$

\begin{figure}[h]
    \centering
    \begin{subfigure}[b]{0.45\textwidth}
        \includegraphics[width=\textwidth]{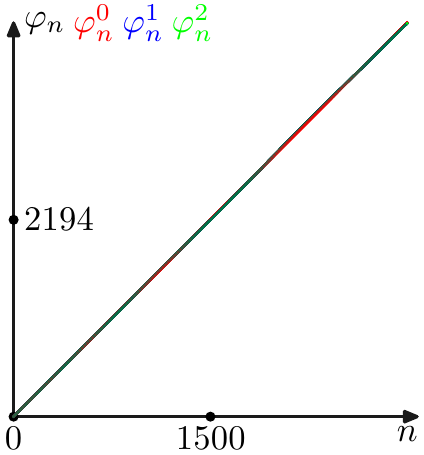}
        \caption{}
    \end{subfigure}
    \begin{subfigure}[b]{0.45\textwidth}
        \includegraphics[width=\textwidth]{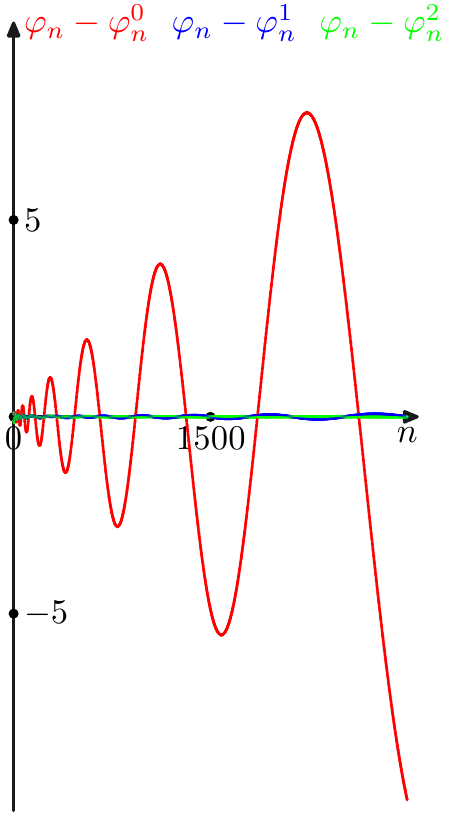}
        \caption{}
    \end{subfigure}
    \caption{For the polynomial $P(z)=0.25z+0.5z^2+0.25z^3$, Tailor coefficients $\vp_n$ (black curve) and their approximations $\vp_n^M$ (red, blue, and green curves) for $M=0,1,2$ are plotted in (a), the differences $\vp_n^M-\vp_n$ (red, blue, and green curves) are plotted in (b).}\lb{fig6}
\end{figure}

{\section{General case}\lb{sec3}} 

In this section, we consider some methods which allow us to determine the complete asymptotics of $\vp_n$ for $n\to+\iy$, up to terms that tend to $0$. However, this section's explanation is rough enough: many assumptions are not mentioned. This Section is just a preparation for further research. While this content is beyond our current needs, it allows us to estimate the possibilities of the above approach. Following the ideas presented in the first part of this article, one may write
\[\lb{300}
 \Psi(z)^{\frac{2\pi\mathbf{i}m+\ln p_1}{\ln E}}=\sum_{j=0}^{+\iy}R_j(2\pi\mathbf{i}m)(1-z)^{\frac{2\pi\mathbf{i}m+\ln p_1}{\ln E}+j},
\]
where the polynomials $R_j(x)$ of degree $j$ have real coefficients determined from the functional equation for $\Psi$, see \er{010}. The next step is the use of expansions
\[\lb{301}
 (1-z)^{\frac{2\pi\mathbf{i}m+\ln p_1}{\ln E}+j}=\sum_{n=0}^{+\iy}(-1)^n\binom{\frac{2\pi\mathbf{i}m+\ln p_1}{\ln E}+j}{n}z^n
\]
in \er{300}, which leads to
\[\lb{302}
 \Psi(z)^{\frac{2\pi\mathbf{i}m+\ln p_1}{\ln E}}=\sum_{n=0}^{+\iy}z^n\sum_{j=0}^{+\iy}R_j(2\pi\mathbf{i}m)(-1)^n\binom{\frac{2\pi\mathbf{i}m+\ln p_1}{\ln E}+j}{n}.
\]
Then, using \er{302} along with \er{018} and \er{006}, we deduce that
\[\lb{303}
 \vp_n=\sum_{j=0}^{+\iy}\sum_{m=-\iy}^{+\iy}\theta_mR_j(2\pi\mathbf{i}m)(-1)^n\binom{\frac{2\pi\mathbf{i}m+\ln p_1}{\ln E}+j}{n}.
\]
The next step is the approximation of the binomial coefficients
\[\lb{304}
 (-1)^n\binom{k}{n}=\frac{n^{-k-1}}{\Gamma(-k)}\sum_{r=0}^{+\iy}\frac{S_{2r}(k)}{n^r},
\]
with polynomials $S_{2r}$ of the degree $2r$. These polynomials can be determined explicitly from the identity
\[\lb{305}
 \lt(1-\frac{k+1}{n+1}\rt)(-1)^n\binom{k}{n}=(-1)^{n+1}\binom{k}{n+1},
\]
which leads to
\begin{multline}\lb{306}
 \lt(1-\frac{k+1}{n}+\frac{k+1}{n^2}-...\rt)\sum_{r=0}^{+\iy}\frac{S_{2r}(k)}{\Gamma(-k)n^{r+k+1}}=\sum_{r=0}^{+\iy}\frac{S_{2r}(k)}{\Gamma(-k)(n+1)^{r+k+1}}=\\
 \sum_{r=0}^{+\iy}\frac{S_{2r}(k)}{\Gamma(-k)n^{r+k+1}}\sum_{p=0}^{+\iy}\frac{1}{n^p}\binom{-r-k-1}{p}.
\end{multline}
Equating the terms with the same power of $n$ and using $S_0\ev1$, we can determine $S_{2r}(k)$, e.g.
\[\lb{307}
 S_2(k)=\frac{k(k+1)}2,\ \ \ S_{4}(k)=\frac{k(k+1)(k+2)(3k+1)}{24}, ....
\]
Using \er{304} in \er{303} leads to 
\[\lb{308}
 \vp_n=\sum_{j=0}^{+\iy}\sum_{m=-\iy}^{+\iy}\theta_mR_j(2\pi\mathbf{i}m)  
 \sum_{r=0}^{+\iy}\frac{S_{2r}(\frac{2\pi\mathbf{i}m+\ln p_1}{\ln E}+j)}{\Gamma(-\frac{2\pi\mathbf{i}m+\ln p_1}{\ln E}-j)n^{r+\frac{2\pi\mathbf{i}m+\ln p_1}{\ln E}+j+1}}.
\]
To obtain the main part of the asymptotics for $n\to+\iy$, we should take the terms with non-negative power of $n$,
\[\lb{309}
 \vp_n\approx\vp^A_n=\sum_{0\le r+j\le\frac{-\ln p_1}{\ln E}-1}n^{\frac{-\ln p_1}{\ln E}-1-r-j}\sum_{m=-\iy}^{m=+\iy}\frac{S_{2r}(\frac{2\pi\mathbf{i}m+\ln p_1}{\ln E}+j)R_j(2\pi\mathbf{i}m)\theta_me^{2\pi\mathbf{i}m\frac{-\ln n}{\ln E}}}{\Gamma(-\frac{2\pi\mathbf{i}m+\ln p_1}{\ln E}-j)}.
\]
Introducing periodic functions
\[\lb{310}
 K_j(z)=\sum_{m=-\iy}^{+\iy}\frac{\theta_me^{2\pi\mathbf{i}mz}}{\Gamma(-\frac{2\pi\mathbf{i}m+\ln p_1}{\ln E}-j)},
\]
we rewrite \er{309} as
\[\lb{311}
 \vp_n\approx\vp^A_n=\sum_{0\le r+j\le\frac{-\ln p_1}{\ln E}-1}n^{\frac{-\ln p_1}{\ln E}-1-r-j}S_{2r}\lt(\frac{\frac{\pa}{\pa x}+\ln p_1}{\ln E}+j\rt)R_j(\frac{\pa}{\pa x})K_j(x)|_{x=\frac{-\ln n}{\ln E}},
\]
since $\pa/\pa x$ leads to the multiplication of the $m$-th Fourier coefficient by $2\pi\mathbf{i}m$. The asymptotic \er{311} contains explicitly the power-terms and periodic factors. As an interesting exercise, one may rewrite \er{310} as a convolution of two functions. In this exercise, it is useful to take into account the symmetry $\theta_m=\ol{\theta_{-m}}$, and the exponential decaying of $\theta_m$, when $m\to\pm\iy$. As it is already mentioned above, \er{311} is very useful in the conceptual context, since it consists of the finite number of explicit terms of the form $n^{\a}{\rm Periodic}(\log_E n)$.

{\section{Conclusion}\lb{sec4}} 

For the solutions of Schr\"oder-type functional equations we provide simple approximations of the corresponding power-series coefficients. A good agreement between the exact values of Taylor coefficients and their approximations is demonstrated by several numerical examples. The complete asymptotic of Taylor coefficients consisting of the finite number of terms $n^{\a}{\rm Periodic}(\log_E n)$ is also discussed.

\section*{Acknowledgements} 
This paper is a contribution to the project M3 of the Collaborative Research Centre TRR 181 "Energy Transfer in Atmosphere and Ocean" funded by the Deutsche Forschungsgemeinschaft (DFG, German Research Foundation) - Projektnummer 274762653. 

\end{document}